\documentclass[draft]{amsart}
\usepackage[latin1]{inputenc}
\usepackage{url,amssymb,bm}

\theoremstyle{plain}
\newtheorem{lemma}{Lemma}
\newtheorem{proposition}[lemma]{Proposition}
\newtheorem{theorem}[lemma]{Theorem}
\newtheorem{corollary}[lemma]{Corollary}

\theoremstyle{definition}

\theoremstyle{remark}
\newtheorem*{remark}{Remark}

\renewcommand{\phi}{\varphi}
\renewcommand{\epsilon}{\varepsilon}
\renewcommand{\rho}{\varrho}
\newcommand{\INumd}{\boldsymbol{\mu}}
\newcommand{\INumdPM}{\INumd^\pm}
\newcommand{\INumdMod}{\tilde{\INumd}^\pm}
\newcommand{\INumdModw}{\tilde{\INumd}^{\pm,w}}
\newcommand{\INumdComb}{\INumd^\circ}
\newcommand{\INumdPerm}{\boldsymbol{\alpha}}
\newcommand{\INumdExt}{\INumdPerm^\circ}
\newcommand{\INumdPermLoc}{\INumdPerm'}
\newcommand{\INht}{\boldsymbol{\kappa}}
\newcommand{\INmt}{\boldsymbol{\tau}}

\newcommand{\txt}[1]{\quad\mbox{#1}\quad}
\newcommand{\kap}{\kappa}

\begin{document}
\title{The UMD constants of the summation operators}

\author{Jörg Wenzel}

\thanks{This article grew out of the author's habilitation thesis,
  which was supported by DFG grant We 1868/1-1.}



\address{Jörg Wenzel, Department of Mathematics and Applied
  Mathematics, University of Pretoria, Pretoria 0002, South Africa}

\email{\url{wenzel@minet.uni-jena.de}}

\date{\today}

\begin{abstract}
  The UMD property of a Banach space is one of the most useful
  properties when one thinks about possible applications. This is in
  particular due to the boundedness of the vector-valued Hilbert
  transform for functions with values in such a space.

  Looking at operators instead of at spaces, it is easy to check that
  the summation operator does not have the UMD property. The actual
  asymptotic behavior however of the UMD constants computed with
  martingales of length $n$ is unknown.

  We explain, why it would be important to know this behavior,
  rephrase the problem of finding these UMD constants and give some
  evidence of how they behave asymptotically.
\end{abstract}

\maketitle

\section{Introduction}
\label{sec:introduction}


A fundamental relation in the theory of Banach spaces is the one
between the Hilbert transform and the unconditionality property for
martingale differences, which was established at the beginning of the
eighties by Burkholder~\cite{bur81b,bur81a} and Bourgain~\cite{bou83}.

To explain this connection, denote the Hilbert transform constant of
an operator $T:X\to Y$ by $\INht(T)$ and its martingale
unconditionality constant by $\INumd(T)$ (see
Sections~\ref{sec:notations} and \ref{sec:umd-property-hilbert} for
precise definitions). Burkholder showed that there is some constant
$c$ such that
\[ \INht(ST)\leq c\INumd(S)\INumd(T),
\]
while it is due to Bourgain that there is some constant
$c$ such that
\[ \INumd(ST)\leq c\INht(S)\INht(T)
\]
for all operators $T:X\to Y$ and $S:Y\to Z$. It is an open problem
even for identity maps of Banach spaces, whether there exists a
constant $c$ such that
\begin{equation}
  \INht(T)/c\leq \INumd(T)\leq c\INht(T),
  \label{eq:11}
\end{equation}
for all linear operators $T:X\to Y$; see
Burkholder~\cite{burkholder:_martin_banac}, especially the problem on
p.~249.

In this paper, we will be interested in the finite summation
operators. For $x=(\xi_k)\in l_1^n$, the \emph{finite summation
  operator} $\Sigma_n:l_1^n\to l_\infty^n$ is defined by
\[ \Sigma_n(x) := \Big(\sum_{k=1}^h\xi_k\Big).
\]
For $\Sigma_{2^n}$, one can easily check that
$\INht(\Sigma_{2^n})\asymp n$ (see
Section~\ref{sec:umd-property-hilbert}), while $\sqrt n\prec
\INumd(\Sigma_{2^n})\prec n$. So if one could show that indeed
$\INumd(\Sigma_{2^n})\asymp \sqrt n$, a relation like~(\ref{eq:11})
could not hold.

From a different viewpoint, the finite summation operators are also
used to characterize superreflexive Banach spaces. Denoting by
$\INumd_n(\Sigma)$ the UMD constant of the infinite summation operator
$\Sigma:l_1\to l_\infty$ computed with martingales of length at most
$n$, if one could show, that $\INumd_n(\Sigma)\asymp n$, this would
establish that every non superreflexive Banach space has
$\INumd_n(X)\asymp n$ and give a nice characterization of
superreflexive Banach spaces.

For both these alternatives, it would be extremely important to know
the martingale unconditionality constants of the finite summation
operators. In this paper, I want to approach this problem, simplify it
and reduce it to a question about a certain matrix norm. We cannot
actually compute these constants, but in the last section, I dare to
make a conjecture based on computer calculations.


\section{Notations}
\label{sec:notations}

For $k=1,2,\dots$, the \emph{dyadic intervals}
\[ \Delta_k^{(i)}:=\Big[\frac {i}{2^k},\frac {i+1}{2^k}\Big)
\txt{where $i=0,\dots,2^k-1$,}
\]
generate the dyadic $\sigma$-algebra denoted by $\mathcal{F}_k$.

For a Banach space $X$, we consider \emph{dyadic martingales}
$(f_1,\dots,f_n)$, defined on $[0,1)$, taking values in $X$, and
adapted to the \emph{dyadic filtration}
$\mathcal{F}_1\subseteq\dots\subseteq\mathcal{F}_n$. We let
$f_0\equiv0$ and denote by $d_k:=f_k-f_{k-1}$ the \emph{differences}
or increments of this martingale.

Given $t\in[0,1)$, we let $\Delta_k(t)$ be the dyadic interval of
length $2^{-k}$ containing $t$, and $\Delta'_k(t)$ its sibling, i.~e.
\[ \Delta'_k(t) := \Delta_{k-1}(t)\setminus\Delta_k(t).
\]

By
\[ \|f|L_p\| := \Big(\int_0^1\|f(t)\|^p\,dt\Big)^{1/p}
\]
we denote the $L_p$-norm of a function $f:[0,1)\to X$.

When dealing with two sequences $(\alpha_n)$ and $(\beta_n)$, we will
use the notations
\[ \alpha_n \prec \beta_n, \quad
\alpha_n\succ \beta_n, \txt{and}
\alpha_n\asymp\beta_n
\]
to indicate that there exists a constant $c>0$ independent of $n$,
such that
\[ \alpha_n \leq c \beta_n, \quad
\alpha_n\geq c\beta_n, \txt{and}
\alpha_n/c\leq\beta_n\leq c\alpha_n,
\]
respectively. In the case $\alpha_n\asymp\beta_n$, we say that the two
sequences are \emph{asymptotically equivalent} or simply
\emph{equivalent} for short.

\section{The UMD property, the Hilbert transform, and
  superreflexivity}
\label{sec:umd-property-hilbert}

For $n\in\mathbb N$ the $n$-th \emph{UMD norm} $\INumd_n(T)$ of an
operator $T:X\to Y$ is the least constant $c\geq0$ such that
\begin{equation}\label{eq:1}
  \Big\| \sum_{k=1}^n \epsilon_k Td_k \Big| L_2 \Big\| \leq c
  \Big\| \sum_{k=1}^n d_k \Big|L_2\Big\|
\end{equation}
for all $X$-valued differences $d_1,\dots,d_n$ of dyadic martingales
and any sequence $\epsilon_1,\dots,\epsilon_n$ of signs.

We let $\INumd(T):=\sup_n\INumd_n(T)$ if this supremum is finite. In
this case, we call $T$ a \emph{UMD operator}.

For $n\in\mathbb{N}$ the $n$-th \emph{Hilbert transform norm}
$\INht_n(T)$ of an operator $T:X\to Y$ is the least constant $c\geq0$
such that
\begin{equation*}
  \Big( \sum_{k=1}^n \Big\| \sum_{\substack{
      h=1\\h\neq k}}^n
  \frac{Tx_h}{h-k}\Big\|^2 \Big)^{1/2} \leq c\Big( \sum_{k=1}^n
  \|x_k\|^2 \Big)^{1/2}
\end{equation*}
for all $n$-tuples of elements $x_1,\dots,x_n\in X$.

We let $\INht(T):=\sup_n \INht_n(T)$ if this supremum is finite. In this
case we call $T$ a \emph{Hilbert transform operator}.

Letting $x_k=e_k$ be the $k$th unit vector in $l_1^{2^n}$, it can
easily be seen that $\INht_{2^n}(\Sigma_{2^n}) \succ n$, while the
reverse estimate $\INht_{2^n}(T)\prec n\|T\|$ is valid for any
operator $T:X\to Y$.

It follows from the estimate $\INht(l_q)\leq cq$ (see
Pichorides~\cite{pich72}) that in fact
\[ \INht(\Sigma_{2^n})\asymp\INht_{2^n}(\Sigma_{2^n})\asymp n.
\]

We now turn our attention to the connection of UMD norms and
super weakly compact (i.~e.~superreflexive) operators.

For $n\in\mathbb N$ the $n$-th \emph{martingale type $2$ norm}
$\INmt_n(T)$ of an operator $T:X\to Y$ is the least constant $c\geq0$
such that
\begin{equation*}
  \Big\| \sum_{k=1}^n  Td_k \Big| L_2 \Big\| \leq c
  \Big( \sum_{k=1}^n \|d_k|L_2\|^2 \Big)^{1/2}
\end{equation*}
for all $X$-valued differences $d_1,\dots,d_n$ of dyadic martingales.

It follows from
\[ \|d_k|L_2\| \leq 2 \Big\|\sum_{k=1}^n d_k \Big| L_2\Big\|
\]
that $\INumd_n(T)\leq 2\sqrt n\INmt_n(T)$. Now every super weakly
compact operator satisfies $\INmt_n(T)/\sqrt n\to 0$ (see
Wenzel~\cite[Thm. 2]{wenzel02:_unifor}), so it follows that every
super weakly compact operator satisfies
\[ \INumd_n(T)/n\to 0.
\]

Since every non superreflexive operator $T$ uniformly factors the
finite summation operators (see James~\cite{jam72c}), we easily obtain
\[ \INumd_n(T) \succ \INumd_n(\Sigma_N)
\]
for all $N$ and all non superreflexive operators $T$. Since every
dyadic martingale of length $n$ actually only takes finitely many
values, we have $\lim_{N\to\infty}\INumd_n(\Sigma_N) =
\INumd(\Sigma)$, where $\Sigma$ denotes the infinite summation
operator on $l_1$.

Summarizing, if $\INumd_n(\Sigma) \asymp n$ then an operator $T$ is
super weakly compact if and only if $\INumd_n(T)/n\to 0$. If on the
other hand $\INumd(\Sigma_{2^n}) \asymp \sqrt n$, then we cannot have
\[ \INumd(T) \leq c \INht(T)
\]
for all linear operators $T$.

\section{Dyadic addition and the function $\kap$}
\label{sec:dyadic-addition}

As it will turn out, a key role in the calculation of
$\INumd_n(\Sigma_{2^n})$ will be played by the dyadic addition and a
certain function, which we will denote by $\kap$. The purpose of this
section is to define these concepts and collect some of their
properties.

We let
\begin{align*}
  0\oplus0=1\oplus1&:=0,\quad&
  1\oplus0=0\oplus1&:=1.
\end{align*}
Given two non negative integers $i$ and $j$ with dyadic expansion
\[ i=\sum_{k=1}^\infty i_k2^{k-1}
\txt{and}
j=\sum_{k=1}^\infty j_k2^{k-1},
\]
where $i_k,j_k\in\{0,1\}$, we let
\[ i\oplus j := \sum_{k=1}^\infty (i_k\oplus j_k) 2^{k-1}.
\]

We denote by $\kap(i)$ the number of binary digits of $i$, that is
\[ \kap(i) :=
\begin{cases}
  \max\{k:i_k\neq0\} & \txt{if $i\neq0$,}\\
  2                  & \txt{if $i    =0$.}
\end{cases}
\]
The reason for the choice $\kap(0)=2$ will become clear in
Lemma~\ref{lem:1}.

We collect here some facts about the function $\kap$.
\begin{proposition}
  If $i,j=1,\dots,2^n-1$ such that $i\neq j$, then
  \begin{flalign}
    \label{eq:3}
    & 2^{\kap(i)-1} \leq i < 2^{\kap(i)} \quad\text{when $i\neq0$,}&&
    \\
    \label{eq:4}
    & \kap(i\oplus j) = \min\{k: (l-1)2^k \leq i,j < l2^k \quad
      \text{for some $l$}\}, \\
    \label{eq:5}
    & i2^{-n} \in \Delta'_k(j2^{-n}) \iff \kap(i\oplus j)=n-k+1, \\
    \label{eq:9}
    & i<j
    \iff j_{\kap(i\oplus j)} = 1
    \iff i_{\kap(i\oplus j)} = 0,
  \end{flalign}
  Concerning the relation of $\kap(i)$ and $\kap(j)$, we have the
  following formulas. To avoid problems with the exceptional case
  $\kap(0)$, we assume here that $k$ is greater than two.
  \begin{flalign}
    \label{eq:6}
    & \kap(i)=k,\ \kap(j)=k \implies \kap(i\oplus j) < k, &&\\
    \label{eq:7}
    & \kap(i)=k,\ \kap(j)<k \implies \kap(i\oplus j) = k, &&\\
    \label{eq:8}
    & \kap(i)<k,\ \kap(j)<k \implies \kap(i\oplus j) < k, &&
  \end{flalign}
\end{proposition}
\begin{proof}
  Inequality~(\ref{eq:3}) is basically the definition of $\kap(i)$.

  Let
  \begin{equation*}\label{eq:13}
    \kappa=\min\{k: (l-1)2^k \leq i,j < l2^k \quad
    \text{for some $l$}\}.
  \end{equation*}
  Then it follows that
  \[ i=(l-1)2^\kappa + \sum_{k=1}^\kappa i_k2^{k-1}
  \txt{and}
  j=(l-1)2^\kappa + \sum_{k=1}^\kappa j_k2^{k-1}.
  \]
  This implies that
  \[ (l-1)2^\kappa = \sum_{k>\kappa} i_k2^{k-1} =  \sum_{k>\kappa}
  j_k2^{k-1}
  \]
  from which we get $i_k=j_k$ for $k>\kappa$. On the other hand, if
  $i_\kappa=j_\kappa$, then
  \[ i=(l-1)2^\kappa + i_\kappa2^{k-1} + \sum_{k=1}^{\kappa-1}
  i_k2^{k-1}
  \txt{and}
  j=(l-1)2^\kappa + i_\kappa2^{k-1} + \sum_{k=1}^{\kappa-1}
  j_k2^{k-1},
  \]
  which implies that for $l'=2l-1+i_\kappa$ we get
  \[ (l'-1)2^{\kappa-1} \leq i,j <
  l'2^{\kappa-1}
  \]
  contradicting the minimality of $\kappa$. Hence $i_\kappa\neq
  j_\kappa$, which means that $\kap(i\oplus j)=\kappa$ and
  proves~(\ref{eq:4}).

  Note that $i2^{-n}\in\Delta'_k(j2^{-n})$ implies
  \[ i2^{-n},j2^{-n} \in \Delta_{k-1}^{(l)}
  \]
  for some $l\in\mathbb{N}$ and no dyadic interval with smaller length
  will contain both $i2^{-n}$ and $j2^{-n}$. This shows~(\ref{eq:5}).

  To see~(\ref{eq:9}), write $\kappa=\kap(i\oplus j)$ and note that
  \[ j-i = \sum_{k=1}^\kappa (j_k-i_k)2^{k-1}
  \begin{cases}
    > 0 & \txt{if $j_k=1$ and $i_k=0$,}\\
    < 0 & \txt{if $j_k=0$ and $i_k=1$.}
  \end{cases}
  \]

  Assume now that $k=\kap(i)=\kap(j)\geq 3$.  Then for $k'>k$ we get
  $i_{k'}=j_{k'}=0$, that is $(i\oplus j)_{k'}=0$ but also $i_k=
  j_k=1$ so that $(i\oplus j)_k = 0$.  That is $\kap(i\oplus j)<k$.
  This shows~(\ref{eq:6}).

  To prove~(\ref{eq:7}) assume that $k=\kap(i)>\kap(j)$. Then for
  $k'>k$ we get $i_{k'}=0$ and since $k'>k>\kap(j)$ also $j_{k'}=0$ so
  that $(i\oplus j)_{k'}=0$. But we also have $i_k=1$ and $j_k=0$ so
  that $\kap(i\oplus j)=k$.

  Formula~(\ref{eq:8}) follows by combining~(\ref{eq:6})
  and~(\ref{eq:7}).
\end{proof}

The following recursive relation is the main reason to let
$\kap(0)=2$.
\begin{lemma}\label{lem:1}
  For all $i=0,1,\dots$ we have
  \[ (-2)^{-\kap(2i\oplus1)} + (-2)^{-\kap(2i)} = -(-2)^{-\kap(i)}.
  \]
\end{lemma}
\begin{proof}
  For $i\neq0$ the relation follows from
  \[ \kap(2i\oplus1) = \kap(2i) = 1+\kap(i).
  \]
  For $i=0$ the assertion is easily checked using the definition and
  this is, where the choice $\kap(0)=2$ plays a role.
\end{proof}

\section{Equivalent UMD norms}
\label{sec:equivalent-umd-norms}

In this section we will define several sequences of ideal norms
related to the unconditionality of martingale differences and prove
their asymptotic equivalence. The final outcome will be an ideal norm
defined with the help of the matrix $\big((-2)^{-\kap(i\oplus
  j)}\big)$.

Admitting only one special sequence of signs in the
definition~(\ref{eq:1}) of the UMD norm, we obtain the following ideal
norm. For $n\in\mathbb N$ the ideal norm $\INumdPM_n(T)$ of an
operator $T:X\to Y$ is the least constant $c\geq0$ such that
\begin{equation*}
  \Big\| \sum_{k=1}^n (-1)^k Td_k \Big| L_2 \Big\| \leq c
  \Big\| \sum_{k=1}^n d_k \Big|L_2\Big\|
\end{equation*}
for all $X$-valued differences $d_1,\dots,d_n$ of dyadic
martingales.

As a further specialization, it is sometimes convenient to consider
the martingale transform
\[ (d_k) \mapsto (2d_{2k-1}-d_{2k}),
\]
which has the advantage that the value of $2d_{2k-1}(t)-d_{2k}(t)$
depends on all the values $f_n(s)$ with
$s\in\Delta_{2k-2}(t)\setminus\Delta_{2k}(t) =
\Delta_{2k-1}'(t)\cup\Delta_{2k}'(t)$, which are disjoint sets for
$k=1,\dots,n$. On the other hand, the value of $d_k(t)$ depends on all
the values $f_n(s)$ with $s\in\Delta_{k-1}(t)$, which are sets
contained in each other.

Therefore, for $n\in\mathbb{N}$ we define $\INumdMod_n(T)$ as the
least constant $c\geq0$ such that
\[ \Big\| \sum_{k=1}^n T(2d_{2k-1} - d_{2k}) \Big| L_2 \Big\| \leq c
\Big\| \sum_{k=1}^{2n} d_k \Big|L_2\Big\|
\]
for all $X$-valued differences $d_1,\dots,d_{2n}$ of dyadic
martingales.

A weaker estimate is obtained by replacing the $L_2$-norm on left by
the $L_1$-norm and on the right by the $L_\infty$-norm. For
$n\in\mathbb{N}$ let $\INumdModw_n(T)$ be the least constant $c\geq0$
such that
\[ \Big\| \sum_{k=1}^n T(2d_{2k-1} - d_{2k}) \Big| L_1 \Big\|
\leq c \Big\| \sum_{k=1}^{2n} d_k \Big|L_\infty\Big\|
\]
for all $X$-valued differences $d_1,\dots,d_{2n}$ of dyadic
martingales.

Finally, for $n\in\mathbb{N}$ let $\INumdComb_n(T)$ be the least
constant $c\geq0$ such that
\[  \frac1{2^n} \sum_{i=0}^{2^n-1} \Big\|
\sum_{j=0}^{2^n-1}
(-2)^{-\kap(i\oplus j)} Tx_j\Big\| \leq c \sup_j\|x_j\|,
\]
for all $x_0,\dots,x_{2^n-1}\in X$.

\begin{lemma}\label{lem:2}
  The sequence $(\INumdComb_n(T))$ is monotonically increasing.
\end{lemma}
\begin{proof}
  Given $x_0,\dots,x_{2^n-1}$, let
  \[ x_{2j}' = x_{2j+1}' := x_j.
  \]
  It follows from Lemma~\ref{lem:1} that
  \[ \frac1{2^n} \sum_{i=0}^{2^n-1} \Big\|
  \sum_{j=0}^{2^n-1}
  (-2)^{-\kap(i\oplus j)} Tx_j\Big\| =
  \frac1{2^{n+1}} \sum_{i=0}^{2^{n+1}-1} \Big\|
  \sum_{j=0}^{2^{n+1}-1}
  (-2)^{-\kap(i\oplus j)} Tx'_j\Big\|
  \]
  from where the monotonicity is immediately clear.
\end{proof}

\begin{theorem}\label{thr:1}
  All of the UMD norms introduced above are asymptotically equivalent.
  We have
  \[ \INumd_n \asymp \INumdPM_n \asymp \INumdMod_n \asymp
  \INumdModw_n \asymp \INumdComb_n.
  \]
\end{theorem}
\begin{proof}
  The equivalence of $\INumd_n$ and $\INumdPM_n$ was proved by the
  author in~\cite{wenzel97:_ideal_umd}.

  To see the equivalence of $\INumdPM_n$ and $\INumdMod_n$, write
  \[ 2d_{2k-1}-d_{2k} = \frac32 (d_{2k-1}-d_{2k}) + \frac12
  (d_{2k-1}+d_{2k}).
  \]
  It follows that
  \[ \Big\| \sum_{k=1}^n T(2d_{2k-1} - d_{2k}) \Big| L_2 \Big\| \leq
  \frac32 \INumdPM_{2n}(T) \Big\| \sum_{k=1}^{2n} d_k \Big|L_2\Big\| +
  \frac12 \|T\|\Big\| \sum_{k=1}^{2n} d_k \Big|L_2\Big\|.
  \]
  This implies
  \[ \INumdMod_n(T)\leq 2\INumdPM_{2n}(T).
  \]
  It can now easily be verified that $\INumdPM_{2n}(T)\leq
  3\INumdPM_n(T)$, see Wenzel~\cite[Prop.~2]{wenzel97:_ideal_umd} and
  consequently
  \[ \INumdMod_n(T)\leq 6\INumdPM_n(T).
  \]

  On the other hand
  \[ d_{2k-1}-d_{2k} = \frac23 (2d_{2k-1}-d_{2k}) - \frac13
  (d_{2k-1}+d_{2k})
  \]
  implies
  \[ \Big\| \sum_{k=1}^{2n} (-1)^k Td_k \Big| L_2 \Big\| \leq
  \frac23 \INumdMod_n(T) \Big\| \sum_{k=1}^{2n} d_k \Big|L_2\Big\| +
  \frac13 \|T\|\Big\| \sum_{k=1}^{2n} d_k \Big|L_2\Big\|.
  \]
  Therefore using the obvious monotonicity of $\INumdPM_n(T)$ we get
  \[ \INumdPM_n(T) \leq \INumdPM_{2n}(T)\leq \INumdMod_n(T),
  \]
  which proves that $\INumdPM_n\asymp\INumdMod_n$.

  That $\INumdModw_n\leq\INumdMod$ follows from the inequalities
  $\|f|L_1\|\leq \|f|L_2\|\leq\|f|L_\infty\|$. The reverse estimate
  can be shown using an extrapolation technique that has its roots in
  Burkholder\slash Gundy~\cite{burkholder70:_extrap} and has been used
  in several places, see Hitczenko~\cite{hitczenko90:_upper_l},
  Geiss~\cite[Theorem~1.7]{geiss97:_bmo}, Pietsch\slash
  Wenzel~\cite[7.2.9]{pietsch98:_orthon_banac}, or
  Wenzel~\cite[Theorem~1, App.~A, p.~58]{wenzel03:_haar_banac}.

  To see the last equivalence, we use the identity
  \begin{equation}
    \sum_{k=1}^n \big(2d_{2k-1}(t)-d_{2k}(t)\big) =
    \sum_{k=1}^{2n} (-2)^k \int_{\Delta'_k(t)} f_{2n}(s)\,ds,
    \label{eq:10}
  \end{equation}
  which follows from the definition of the conditional expectation and
  makes the use of the differences $2d_{2k-1}-d_{2k}$ so useful.
  Denoting by $x_j$ the constant value of $f_{2n}$ on the interval
  $\Delta_{2n}^{(j)}$ for $j=0,\dots,2^{2n}-1$ it follows that for
  $t\in[0,1)$ we have
  \[ \int_{\Delta'_k(t)} f_{2n}(s)\, ds = \frac1{2^{2n}}
  \sum_{j\in\mathbb{N}_k(t)} x_j,
  \]
  where
  \[ \mathbb{N}_k(t) :=
  \big\{ j=0,\dots,2^{2n}-1:\Delta_{2n}^{(j)}\subseteq\Delta'_k(t)
  \big\}.
  \]
  Now from~(\ref{eq:10}) it follows that
  \[ \sum_{k=1}^n T\big(2d_{2k-1}(t) - d_{2k}(t)\big) =
  \sum_{k=1}^{2n} (-2)^{k-2n} \sum_{j\in\mathbb{N}_k(t)} Tx_j.
  \]
  Since $\mathbb{N}_k(t)=\mathbb{N}_k(\frac i{2^{2n}})$ for
  $t\in\Delta_{2n}^{(i+1)}$, we obtain
  \begin{equation*}
    \Big\|\sum_{k=1}^nT\big(2d_{2k-1} - d_{2k}\big)\Big|L_1\Big\| =
    \frac1{2^{2n}} \sum_{i=0}^{2^{2n}-1} \Big\|
    \sum_{k=1}^{2n} (-2)^{k-2n} \sum_{j\in\mathbb{N}_k(\frac i{2^{2n}})}
    Tx_j \Big\|.
  \end{equation*}
  By definition of the sets $\mathbb{N}_k(t)$ and using~(\ref{eq:5})
  we obtain
  \[ j\in\mathbb{N}_k\Big(\frac i{2^{2n}}\Big) \iff \frac j{2^{2n}}
  \in \Delta'_k\Big(\frac i{2^{2n}}\Big) \iff \kap(i\oplus j)=2n-k+1
  \txt{and $i\neq j$.}
  \]
  So we can continue as
  \begin{align*}
    \Big\|\sum_{k=1}^nT\big(2d_{2k-1} - d_{2k}\big)\Big|L_1\Big\|
    &=
    \frac1{2^{2n}} \sum_{i=0}^{2^{2n}-1} \Big\| \sum_{k=1}^{2n}
    (-2)^{k-2n} \sum_{j:\kap(i\oplus j)=2n-k+1} Tx_j \Big\| \\
    &=
    \frac2{2^{2n}} \sum_{i=0}^{2^{2n}-1} \Big\| \sum_{k=1}^{2n}
    \sum_{\substack{j:\kap(i\oplus j)=k\\j\neq i}} (-2)^{-\kap(i\oplus
    j)}Tx_j \Big\|
  \end{align*}
  Finally it is clear that
  \[ \bigcup_{k=1}^{2n} \{j:\kap(i\oplus j)=k\} =
  \{0,\dots,2^{2n}-1\}\setminus\{i\},
  \]
  so that
  \[ \Big\|\sum_{k=1}^nT\big(2d_{2k-1} - d_{2k}\big)\Big|L_1\Big\| =
  \frac2{2^{2n}} \sum_{i=0}^{2^{2n}-1} \Big\|
  \sum_{\substack{j=0\\j\neq i}}^{2^{2n}-1} (-2)^{-\kap(i\oplus
    j)}Tx_j \Big\|.
  \]
  On the other hand
  \[ \Big\| \sum_{k=1}^{2n} d_k \Big|L_\infty\Big\| =
  \|f_{2n}|L_\infty\| = \sup_j\|x_j\|.
  \]
  Since
  \[ \Bigg| \Big\|
  \sum_{\substack{j=0\\j\neq i}}^{2^{2n}-1} (-2)^{-\kap(i\oplus
    j)}Tx_j \Big\| -
  \Big\|
  \sum_{j=0}^{2^{2n}-1} (-2)^{-\kap(i\oplus
    j)}Tx_j \Big\| \Bigg| \leq \frac14 \|T\| \sup_j\|x_j\|,
  \]
  these facts imply that
  \[ \frac49 \INumdModw_n(T) \leq \INumdComb_{2n}(T) \leq \frac34
  \INumdModw_n(T).
  \]
  Now using again that $\INumdModw_n\asymp\INumdPM_n \asymp
  \INumdPM_{2n}\asymp \INumdModw_{2n}$, the monotonicity of
  $\INumd_n\asymp\INumdModw_n$ and the monotonicity of $\INumdComb_n$
  (Lemma~\ref{lem:2}), we obtain the complete equivalence.
\end{proof}

The next theorem specializes the UMD norms to the case of the finite
summation operators. To do so, we introduce two further sequences of
numbers.

For $n\in\mathbb{N}$ let
\begin{equation*}
  \INumdExt_n:=
  \sup_{\pi,(\epsilon_j)}
  \frac1{2^n} \sum_{i=0}^{2^n-1} \sup_{0\leq h<2^n} \Big|
  \sum_{j:\pi(j)\leq h}
  (-2)^{-\kap(i\oplus j)} \epsilon_j \Big|.
\end{equation*}
where the supremum is taken over all maps
$\pi:\{0,\dots,2^n-1\}\to\{0,\dots,2^n-1\}$ and all $\epsilon_j=\pm1$.

For $n\in\mathbb{N}$ let
\begin{equation*}
  \INumdPerm_n:=
  \sup_{\pi}
  \frac1{2^n} \sum_{i=0}^{2^n-1} \sup_{0\leq h<2^n} \Big|
  \sum_{j:\pi(j)\leq h}
  (-2)^{-\kap(i\oplus j)} \Big|.
\end{equation*}
where the supremum is taken over all permutations
$\pi$ of the set $\{0,\dots,2^n-1\}$.

\begin{theorem}\label{thr:2}
  The quantities introduced above are asymptotically equivalent to the
  UMD norm of the summation operators $\Sigma_{2^n}$. We have
  \[ \INumd_n(\Sigma_{2^n}) \asymp \INumdExt_n \asymp \INumdPerm_n.
  \]
\end{theorem}
\begin{proof}
  In the case of an operator $T$ starting in $l_1^{2^n}$, by an
  extreme point argument the vectors $x_j\in l_1^{2^n}$ appearing in
  the definition of $\INumdComb_n(T)$ can be taken as signed unit
  vectors in $l_1^{2^n}$. That is, there exists a map
  $\pi:\{0,\dots,2^n-1\}\to\{0,\dots,2^n-1\}$ and signs
  $\epsilon_j=\pm1$, such that
  \[ \INumdComb(\Sigma_{2^n}) = \frac1{2^n} \sum_{i=0}^{2^n-1} \Big\|
  \sum_{j=0}^{2^n-1} \Sigma_{2^n}\epsilon_j e_{\pi(j)} \Big|
  l_\infty^{2^n} \Big\|.
  \]
  This shows that
  \[ \INumdComb_n(\Sigma_{2^n}) = \INumdExt_n.
  \]

  To see the equivalence of $\INumdExt_n$ and $\INumdPerm_n$,
  i.~e.~that we can actually assume that $\pi$ is a permutation and
  $\epsilon_i=1$, we prepare some lemmas.
\end{proof}

\begin{lemma}\label{lem:3}
  Given a map $\pi:\{0,\dots,2^n-1\}\to\{0,\dots,2^n-1\}$ there is a
  permutation $\rho:\{0,\dots,2^n-1\}\to\{0,\dots,2^n-1\}$ such that
  for any function $f:\{0,\dots,2^n-1\}\to\mathbb{R}$ we have
  \[ \sup_{0\leq h<2^n}
  \Big| \sum_{j:\pi(j) \leq h} f(j) \Big|
  \leq
  \sup_{0\leq h<2^n}
  \Big| \sum_{j:\rho(j)\leq h} f(j) \Big|.
  \]
\end{lemma}
\begin{proof}
  For $h=0,\dots,2^n-1$ let $\mathbb{F}_h:=\pi^{-1}(h) =
  \{j:\pi(j)=h\}$.

  For $j\in\mathbb{F}_h$ define $\rho(j)$ by
  \[ \rho(j) := |\mathbb{F}_0| + \dots + |\mathbb{F}_{h-1}| +
  |\{j'\in\mathbb{F}_h: j' < j\}|.
  \]
  That is, $\rho(j)$ is obtained by counting all the indices that have
  smaller images than $j$ under $\pi$ plus all the indices that have
  the same image under $\pi$ and are smaller than $j$. Note that
  $\rho$ is injective, hence a permutation.

  For every $j,h\in\{0,\dots,2^n-1\}$ we have
  \begin{align*}
    \pi(j)\leq h
    \iff&
    j\in\mathbb{F}_0\cup\dots\cup\mathbb{F}_h\\
    \iff&
    \rho(j) < |\mathbb{F}_0|+\dots+|\mathbb{F}_h|
    \iff
    \rho(j)\leq|\mathbb{F}_0|+\dots+|\mathbb{F}_h|-1.
  \end{align*}
  So given $h$, we either have
  $|\mathbb{F}_0|+\dots+|\mathbb{F}_h|=0$, in which case
  $\{j:\pi(j)\leq h\}$ is empty and
  \[
  \Big| \sum_{j:\pi(j) \leq h } f(j) \Big|
  =
  0\leq
  \sup_{0\leq h<2^n}
  \Big| \sum_{j:\rho(j)\leq h} f(j) \Big|,
  \]
  or we have $h' := |\mathbb{F}_0|+\dots+|\mathbb{F}_h|-1 \in
  \{0,\dots,2^n-1\}$, in which case
  \[ \{j:\pi(j)\leq h\} = \{j:\rho(j)\leq h'\}
  \]
  and also
  \[
  \Big| \sum_{j:\pi(j) \leq h } f(j) \Big|
  =
  \Big| \sum_{j:\rho(j)\leq h'} f(j) \Big|
  \leq
  \sup_{0\leq h<2^n}
  \Big| \sum_{j:\rho(j)\leq h} f(j) \Big|.
  \]
  Taking the supremum over all $h$ on the left hand side proves the
  assertion.
\end{proof}
\begin{lemma}\label{lem:4}
  Given a subset $A\subseteq\{0,\dots,2^n-1\}$ and a permutation $\pi$
  of the set $\{0,\dots,2^n-1\}$, there exists a permutation $\rho$ of
  the same set such that for any function $f:\{0,\dots,2^n-1\} \to
  \mathbb{R}$, we have
  \[ \sup_{0\leq h<2^n}
  \Big| \sum_{\substack{j:\pi(j)\leq h\\j\in A}} f(j) \Big|
  \leq
  \sup_{0\leq h<2^n}
  \Big| \sum_{j:\rho(j)\leq h}                   f(j) \Big|.
  \]
\end{lemma}
\begin{proof}
  Given $h$ we define $h'$ by $\pi(h')=\max\{\pi(j'):j'\in A,
  \pi(j')\leq h\}$.  For $j\in A$ apparently
  \[ \pi(j)\leq h \iff \pi(j)\leq \pi(h').
  \]
  That means, that we can replace the supremum over all $h$ by the
  supremum over $\pi(h')$ with $h'\in A$.

  We now define the permutation $\rho$ by
  \[ \rho(j) =
  \begin{cases}
        |\{j'    \in A:\pi(j')<\pi(j)\}| & \mbox{if $j    \in A$,} \\
    |A|+|\{j'\not\in A:\pi(j')<\pi(j)\}| & \mbox{if $j\not\in A$.}
  \end{cases}
  \]
  That is, if one considers a permutation of $\{0,\dots,2^n-1\}$ as a
  list of the numbers $0,\dots,2^n-1$, to get $\rho$ we first list all
  the numbers of $A$ in the order they appear in the list for $\pi$
  and then all the remaining numbers also in the order they appear in
  the list for $\pi$.

  This permutation preserves monotonicity on $A$, i.~e. for $h'\in A$
  we have
  \[ \big( j\in A \txt{and} \pi(j)\leq \pi(h') \big)
  \Longleftrightarrow \rho(j)\leq\rho(h'),
  \]
  which implies for $h'\in A$ that
  \[ \{j:\pi(j)\leq \pi(h'),j\in A\} = \{j:\rho(j)\leq\rho(h')\}.
  \]
  Therefore
  \begin{equation*}
    \sup_{h'\in A} \Big| \sum_{\substack{
        j:\pi(j)\leq \pi(h')\\j\in A}} f(j) \Big|
    =
    \sup_{h'\in A} \Big| \sum_{
        j:\rho(j)\leq \rho(h')} f(j) \Big|
    \leq
    \sup_{0\leq h<2^n} \Big| \sum_{
        j:\rho(j)\leq h} f(j) \Big|.\qedhere
  \end{equation*}
\end{proof}
We can now finish the proof of Theorem~\ref{thr:2}.
\begin{proof}[Proof of Theorem~\ref{thr:2}. (cont.)]
  We trivially have $\INumdExt_n\geq\INumdPerm_n$. On the other hand,
  given a map $\pi:\{0,\dots,2^n-1\}\to\{0,\dots,2^n-1\}$ and signs
  $\epsilon_j=\pm1$, we first find a permutation $\rho$ according to
  Lemma~\ref{lem:3} such that for all $i$ we have
  \[ \sup_{0\leq h<2^n} \Big| \sum_{j:\pi(j)\leq h}
  (-2)^{-\kap(i\oplus j)} \epsilon_j \Big| \leq
  \sup_{0\leq h<2^n} \Big| \sum_{j:\rho(j)\leq h}
  (-2)^{-\kap(i\oplus j)} \epsilon_j \Big|.
  \]
  Next we let $A_\pm:=\{j:\epsilon_j=\pm1\}$ and obtain permutations
  $\rho_\pm$ according to Lemma~\ref{lem:4} such that for all $i$ we
  have
  \[ \sup_{0\leq h<2^n} \Big| \sum_{\substack{j:\rho(j)\leq h\\j\in
      A_\pm}} (-2)^{-\kap(i\oplus j)} \Big|
  \leq
  \sup_{0\leq h<2^n} \Big| \sum_{j:\rho_\pm(j)\leq h}
  (-2)^{-\kap(i\oplus j)} \Big|,
  \]
  which by
  \[ \Big| \sum_{j:\rho(j)\leq h}
  (-2)^{-\kap(i\oplus j)} \epsilon_j \Big| \leq
  \Big| \sum_{\substack{j:\rho(j)\leq h\\j\in A_+}}
  (-2)^{-\kap(i\oplus j)} \Big| +
  \Big| \sum_{\substack{j:\rho(j)\leq h\\j\in A_-}}
  (-2)^{-\kap(i\oplus j)} \Big|
  \]
  implies that $\INumdExt_n\leq 2 \INumdPerm_n$.
\end{proof}

\section{Special permutations}
\label{sec:special-permutations}

With Theorem~\ref{thr:2} the problem of the computation of the UMD
norm of the summation operators $\Sigma_{2^n}$ is reduced to the
maximization of a certain expression over all possible permutations of
the set $\{0,\dots,2^n-1\}$. For $n\in\mathbb{N}$ and a permutation
$\pi$ of the set $\{0,\dots,2^n-1\}$ let
\begin{equation}\label{eq:12}
  \alpha_n(\pi) := \frac1{2^n} \sum_{i=0}^{2^n-1} \sup_{0\leq h<2^n}
  \Big| \sum_{j:\pi(j)\leq h}
  (-2)^{-\kap(i\oplus j)} \Big|.
\end{equation}
Apparently we get
\[ \INumd_n(\Sigma_{2^n}) \asymp \INumdPerm_n = \sup_\pi \alpha_n(\pi).
\]
In this section, we will take a closer look at the numbers
$\alpha_n(\pi)$ for various permutations $\pi$.

To get a further handle on the numbers $\alpha_n(\pi)$ we will first
analyze the expression
\[ \sum_{j:\pi(j)\leq h} (-2)^{-\kap(i\oplus j)}
= \sum_{\substack{j:\pi(j)\leq h\\j=i}} \frac14 +
\sum_{k=1}^n (-2)^{-k} |\{j:\pi(j)\leq h,j\neq i,\kap(i\oplus j)=k\}|.
\]
For $k\geq 3$, the last sets can be split up further as follows
\begin{multline*}
  \{j:\pi(j)\leq h,j\neq i,\kap(i\oplus j)=k\}
  =\{j:\pi(j)\leq h,\kap(i\oplus j)=k\} = \\
  = \bigcup_{l=1}^n \bigg\{j:
  \pi(j)<h,
  \begin{gathered}
    \kap(i\oplus j)=k,\\
    \kap(h\oplus\pi(j))=l
  \end{gathered}
  \bigg\}
  \cup \{j:\pi(j)=h, \kap(i\oplus j)=k\}.
\end{multline*}
By~(\ref{eq:9}) we have for $\kap(h\oplus\pi(j))=l\geq 3$ that
$\pi(j)<h$ if and only if $h_l=1$, so
\[ \bigg|\bigg\{j:
  \pi(j)<h,
\begin{gathered}
  \kap(i\oplus j)=k,\\
  \kap(h\oplus\pi(j))=l
\end{gathered}
\bigg\}\bigg|
= h_l\cdot \bigg|\bigg\{j:
\begin{gathered}
  \kap(i\oplus j)=k, \\
  \kap(h\oplus\pi(j))=l
\end{gathered}
\bigg\}\bigg|
\]
Moreover, since the exceptional sets for $h\leq 2$ and $l\leq 2$ are
all bounded in size by four and since
\[ \Big|\sum_{k=1}^n (-2)^{-k}\cdot4\Big| \leq 4,
\]
we can write
\begin{equation}
  \alpha_n(\pi) \asymp \frac1{2^n} \sum_{i=0}^{2^n-1} \sup_{0\leq
    h<2^n} \Big| \sum_{k,l=3}^n (-2)^{-k} h_l |A_{kl}^\pi(i,h)| \Big|,
  \label{eq:14}
\end{equation}
where we use the notation
\[ A_{kl}^\pi(i,h) := \{j:\kap(i\oplus j)=k, \kap(h\oplus\pi(j))=l\}.
\]
The following lemma gives the sizes of a simpler version of those
sets.
\begin{lemma}\label{lem:5}
  For any $i=0,\dots2^n-1$ and $k=3,\dots,n$, we have
  \[ \big|\{j:\kap(i\oplus j)=k\}\big|= 2^{k-1}.
  \]
\end{lemma}
\begin{proof}
  Write $j=\sum_{l=1}^nj_l2^{l-1}$. And note that
  \begin{multline*}
    \big|\{j:\kap(i\oplus j)=k\}\big|
    = \big|\{i\oplus j:\kap(i\oplus j)=k\}\big|
    = \big|\{j:\kap(j)=k\}\big| \\
    = \big|\{j:  j_k\neq 0,j_{k+1}=\dots=j_n=0\}\big|
    = \big|\{2^{k-1},\dots,2^k-1\}\big| = 2^{k-1}.\qedhere
  \end{multline*}
\end{proof}

Let us next do the most obvious thing and use the identity permutation
$\pi=\iota$ in the calculation of $\alpha_n(\pi)$. For simplicity we
will write $A_{kl}(i,h):=A_{kl}^\iota(i,h)$ and first determine the size
of these sets.

\begin{lemma}\label{lem:6}
  If $i\neq h$ and $k,l=3,\dots,n$ then
  \[ |A_{kl}(i,h)| =
  \begin{cases}
    2^{k-1} & \mbox{if $l=\kap(i\oplus h)>k$,}\\
    2^{l-1} & \mbox{if $k=\kap(i\oplus h)>l$,}\\
    2^{l-1} & \mbox{if $k=l>\kap(i\oplus h)$,}\\
    0 & \mbox{otherwise.}
  \end{cases}
  \]
\end{lemma}
\begin{proof}
  We distinguish three cases.
  \begin{description}
  \item[Case 1 $k<l$] It follows from~(\ref{eq:7}) that for
    $\kap(i\oplus h)\neq l$ we have $A_{kl}(i,h)=\emptyset$. If on the
    other hand $\kap(i\oplus h)=l$, then again by~(\ref{eq:7}) we have
    that $\kap(i\oplus j)=k$ already implies $\kap(h\oplus j)=l$.
    Therefore $A_{kl}(i,h)=\{j:\kap(i\oplus j)=k\}$. Summarizing, we
    get
    \[ A_{kl}(i,h)=
    \begin{cases}
      \emptyset               & \mbox{if $\kap(i\oplus h)\neq l$,}\\
      \{j:\kap(i\oplus j)=k\} & \mbox{if $\kap(i\oplus h)=l$.}
    \end{cases}
    \]
    The assertion in this case now follows from Lemma~\ref{lem:5}.
  \item[Case 2 $k>l$] We obtain similarly
    \[ A_{kl}(i,h)=
    \begin{cases}
      \emptyset               & \mbox{if $\kap(i\oplus h)\neq k$,}\\
      \{j:\kap(i\oplus j)=l\} & \mbox{if $\kap(i\oplus h)=k$,}
    \end{cases}
    \]
    and the assertion follows again from Lemma~\ref{lem:5}.
  \item[Case 3 $k=l$] Using~(\ref{eq:6}) instead of~(\ref{eq:7}), we
    get in this case
    \[ A_{kk}(i,h)=
    \begin{cases}
      \emptyset               & \mbox{if $\kap(i\oplus h) \geq k$,}\\
      \{j:\kap(i\oplus j)=k\} & \mbox{if $\kap(i\oplus h)<     k$,}
    \end{cases}
    \]
    and the assertion follows once again from Lemma~\ref{lem:5}.\qedhere
  \end{description}
\end{proof}

We can now estimate the quantities $\alpha_n(\iota)$.
\begin{theorem}
  We have $\alpha_n(\iota) \asymp \sqrt n$. Moreover, choosing $h=i$
  in~(\ref{eq:14}) asymptotically maximizes the expression for
  $\alpha_n(\iota)$.
\end{theorem}
\begin{proof}
  Writing
  \begin{align*}
    S_< & := \sum_{k<l} (-2)^{-k} h_l |A_{kl}(i,h)|, \\
    S_> & := \sum_{k>l} (-2)^{-k} h_l |A_{kl}(i,h)|, \\
    S_= & := \sum_{k=l} (-2)^{-k} h_l |A_{kl}(i,h)|
  \end{align*}
  we can split the summation over $k$ and $l$ in~(\ref{eq:14}) into
  three parts. Hence
  \[ \alpha_n(\iota) \asymp \frac1{2^n} \sum_{i=0}^{2^n-1}
  \sup_{0\leq h<2^n} \big(|S_<|+|S_>|+|S_=|\big).
  \]

  It follows from Lemma~\ref{lem:6} that
  \begin{align*}
    S_< & = \sum_{k<\kap(i\oplus h)} (-2)^{-k} h_{\kap(i\oplus h)}
    2^{k-1}, \\
    S_> & = \sum_{\kap(i\oplus h)>l} (-2)^{-\kap(i\oplus h)} h_l
    2^{l-1}, \\
    S_= & = \sum_{k>\kap(i\oplus h)} (-2)^{-k} h_k 2^{k-1}.
  \end{align*}
  The absolute values of $S_<$ and $S_>$ can easily be estimated:
  \[ |S_<|=\Big| h_{\kap(i\oplus h)} \sum_{k<\kap(i\oplus h)}
  (-2)^{-k} 2^{k-1} \Big| \leq \frac12
  \]
  and
  \[ |S_>|\leq 2^{-\kap(i\oplus h)} \sum_{l<\kap(i\oplus h)} 2^{l-1}
  \leq \frac12.
  \]
  Furthermore, if $k>\kap(i\oplus h)$ then $i_k=h_k$. This implies
  \[ S_= = \sum_{k>\kap(i\oplus h)} (-1)^k i_k
  \]
  and the substitution $m=\kap(i\oplus h)$ gives
  \[ \alpha_n(\iota) \asymp \frac1{2^n} \sum_{i=0}^{2^n-1}
  \sup_{1\leq m\leq n} \Big| \sum_{k>m} (-1)^k i_k \Big|.
  \]

  Now, writing $t=i2^{-n}$, we get $i_k=(r_k(t)+1)/2$, where $r_k$
  denotes the $k$th Rademacher function and the sum becomes an
  integral over $t$:
  \[ \alpha_n(\iota) \asymp \int_0^1 \sup_m \Big| \sum_{k>m}
  (-1)^k (r_k(t)+1)\Big|\, dt
  \asymp \int_0^1 \sup_m \Big| \sum_{k>m} (-1)^k r_k(t)\Big|\,
  dt.
  \]
  By Hölder's and Doob's inequality the last integral is bounded by
  \[ \Bigg( \int_0^1 \sup_m \Big| \sum_{k=m+1}^n (-1)^k
  r_k(t) \Big|^2 \,dt \Bigg)^{1/2}
  \leq
  2 \Bigg( \int_0^1 \Big|
  \sum_{k=1}^n (-1)^k r_k(t) \Big|^2 \,dt \Bigg)^{1/2} = 2
  \sqrt{n}.
  \]
  On the other hand, taking $m=1$ and using Khintchine's inequality,
  we obtain
  \[ \alpha_n(\iota) \succ
  \int_0^1  \Big| \sum_{k=2}^n (-1)^k
  r_k(t) \Big|\,dt \succ \sqrt n.
  \]
  This proves the first part of the theorem.

  The moreover part follows by realizing that indeed
  \[ \frac1{2^n} \sum_{i=0}^{2^n-1} \Big| \sum_{k,l=3}^n (-2)^{-k} h_l
  |A_{kl}(i,i)| \Big| \asymp \sqrt n.\qedhere
  \]
\end{proof}

Apparently it is the regularity in the size of the sets
$A_{kl}^\pi(i,h)$, which makes the proof of this theorem work for
$\pi=\iota$.

Another class of permutations for which we can describe the size of
the sets $A_{kl}^\pi(i,h)$ are so called dyadically linear
permutations. Since the method we are going to develop actually works
for a slightly more general class, we will first describe this class
of permutations.

A permutation $\pi$ of $\{0,\dots,2^n-1\}$ is called
\emph{(dyadically) linear}, if
\[ \pi(i\oplus j)=\pi(i)\oplus\pi(j)
\]
for all $i,j\in\{0,\dots,2^n-1\}$. Linear permutations are considered
in \cite[p.~16]{schipp90:_walsh_hilger} under the name
$\mathbb{Z}_2$-linear permutations in order to study relatives of the
Walsh system of functions.

Let us call a permutation $\pi$ of $\{0,\dots,2^n-1\}$ \emph{pseudo
  linear}, if
\begin{equation}
  \kap(\pi(i\oplus j)\oplus\pi(0))=\kap(\pi(i)\oplus\pi(j))
\end{equation}
for all $i,j\in\{0,\dots,2^n-1\}$.

Of course, every linear permutation is pseudo linear, but not
conversely as is seen by the permutation
$\pi:(0,1,2,3,4,5,6,7)\mapsto(0,1,2,3,4,5,7,6)$. Note that for this
permutation $\pi(4)\oplus\pi(6)=4\oplus7=3$, while
$\pi(4\oplus6)=\pi(2)=2$. So we have indeed a bigger class of
permutations. However, this is no longer a group, since the
composition of two pseudo linear permutations need not be pseudo
linear. This can be seen by composing the permutation
$(0,1,2,3,4,5,7,6)$ which is pseudo linear and $(0,4,2,6,1,5,3,7)$,
which is even linear. Also the inverse of a pseudo linear permutation
need not be pseudo linear. An example is the permutation
$(0,6,2,5,3,4,1,7)$.

For pseudo linear permutations the size of the sets $A_{kl}^\pi(i,h)$
behaves quite regular, as we will prove in Proposition~\ref{prop:1}.

In order to prepare the proof of this fact, we define the following
relatives of the sets $A_{kl}^\pi(i,h)$. Let both $\vartriangleleft$
and $\blacktriangleleft$ denote one of the relations $<$ or $=$. We
let
\[ A_{\vartriangleleft k,\blacktriangleleft l}^\pi(i,h)
:=
\{j:
  \kap(i\oplus j)\vartriangleleft n,
  \kap(h\oplus\pi(j))\blacktriangleleft l
\}.
\]
We simply write
\[ A_{\vartriangleleft k,\blacktriangleleft l}^\pi:=
A_{\vartriangleleft k,\blacktriangleleft l}^\pi(0,\pi(0)).
\]
The previously used sets $A_{kl}^\pi(i,h)$ now appear as
$A_{=k,=l}^\pi(i,h)$ for $k,l\geq 3$. We want to replace them in the
estimate of $\alpha_n(\pi)$ by sets of this form that can be handled
easier. The first such replacement works for any permutation.
\begin{lemma}
\[  \alpha_n(\pi) \asymp \frac1{2^n} \sum_{i=0}^{2^n-1} \sup_{0\leq
  h<2^n} \Big| \sum_{k,l=3}^n (-2)^{-k} h_l |A_{<k,=l}^\pi(i,h)|
  \Big|.
\]
\end{lemma}
\begin{proof}
  Since for $k,l= 3,\dots,n$
  \[ A_{kl}^\pi(i,h)=A_{=k,=l}^\pi(i,h)= A_{<k+1,=l}^\pi(i,h)\setminus
  A_{<k,=l}^\pi(i,h)
  \]
  we can write
  \begin{align}\label{eq:15}
    \sum_{k=3}^n (-2)^{-k} |A_{kl}^\pi(i,h)|
    &=
    \sum_{k=3}^n (-2)^{-k} |A_{<k+1,=l}^\pi(i,h)| -
    \sum_{k=3}^n (-2)^{-k} |A_{<k,  =l}^\pi(i,h)|
  \end{align}
  For the first summand we get by an index shift
  \begin{multline*}
    \sum_{k=3}^n (-2)^{-k} |A_{<k+1,=l}^\pi(i,h)| =
    (-2)\sum_{k=3}^n (-2)^{-k} |A_{<k,=l}^\pi(i,h)| + \\
    (-2)^{-n} |A_{<n+1,=l}^\pi(i,h)| - (-2)^{-2} |A_{<3,=l}^\pi(i,h)|.
  \end{multline*}
  But for the last two sets we have the trivial estimate
  \begin{align*}
    |A_{<n+1,=l}^\pi(i,h)| &\leq
    |\{0,1,\dots,2^n-1\}| = 2^n, \\
    |A_{<3,=l}^\pi(i,h)| &\leq
    |\{j:\kap(i\oplus j)<3\}| = 4
  \end{align*}
  so that
  \[ \sum_{k=3}^n (-2)^{-k} |A_{<k+1,=l}^\pi(i,h)| -
  \sum_{k=3}^n (-2)^{-k} |A_{<k,  =l}^\pi(i,h)|
  \asymp -3 \sum_{k=3}^n (-2)^{-k} |A_{<k,=l}^\pi(i,h)|.
  \]
  Multiplication by $h_l$ and summation over $l$ in~(\ref{eq:15}) then
  gives
  \begin{align*}\label{eq:15}
    \sum_{k,l=3}^n (-2)^{-k} h_l |A_{kl}^\pi(i,h)|
    &\asymp
    -3 \sum_{k,l=3}^n (-2)^{-k} h_l |A_{<k,=l}^\pi(i,h)|,
  \end{align*}
  which completes the proof.
\end{proof}

To formulate the next result, we will use the notation
\[ A\oplus B=\{a\oplus b: a\in A, b\in B\}.
\]
If either $A$ or $B$ are empty, we simply let $A\oplus B=\emptyset$.

For pseudo linear permutations we can then prove the following result
about the sets $A_{\vartriangleleft k,\blacktriangleleft l}^\pi(i,h)$.
\begin{lemma}\label{lem:8}
  Let $\pi$ be a pseudo linear permutation, $i,i',h,h'=0,\dots,2^n-1$,
  and $k,l=3,\dots,n$.
  \begin{align*}
    A_{<k,=l}^\pi(i,\pi(h))\oplus A_{<k,=l}^\pi(i',\pi(h'))
    &\subseteq A_{<k,<l}^\pi(i\oplus i',\pi(h\oplus h')), \\
    A_{<k,<l}^\pi(i,\pi(h))\oplus A_{<k,=l}^\pi(i',\pi(h'))
    &\subseteq A_{<k,=l}^\pi(i\oplus i',\pi(h\oplus h')), \\
    A_{=k,<l}^\pi(i,\pi(h))\oplus A_{=k,<l}^\pi(i',\pi(h'))
    &\subseteq A_{<k,<l}^\pi(i\oplus i',\pi(h\oplus h')), \\
    A_{<k,<l}^\pi(i,\pi(h))\oplus A_{=k,<l}^\pi(i',\pi(h'))
    &\subseteq A_{=k,<l}^\pi(i\oplus i',\pi(h\oplus h')).
  \end{align*}
\end{lemma}
\begin{proof}
  Let $j\in A_{<k,=l}^\pi(i,\pi(h))$ and $j'\in
  A_{<k,=l}^\pi(i',\pi(h'))$, i.~e.
  \[ \kap(i\oplus j),\kap(i'\oplus j')<k \txt{and}
  \kap(\pi(h)\oplus\pi(j))=\kap(\pi(h')\oplus\pi(j'))=l.
  \]
  By pseudo linearity we get
  \[ \kap(\pi(h)\oplus\pi(j))=\kap(\pi(h\oplus j)\oplus\pi(0))
  \txt{and}
  \kap(\pi(h')\oplus\pi(j'))=\kap(\pi(h'\oplus j')\oplus\pi(0))
  \]
  So that from relation~(\ref{eq:8}) it follows that $\kap(i\oplus
  i'\oplus j\oplus j')<k$ while from~(\ref{eq:6}) we get
  $\kap(\pi(h\oplus j)\oplus \pi(h'\oplus j'))<l$. Now again using the
  pseudo linearity of $\pi$ twice gives
  \[ \kap(\pi(h\oplus j)\oplus \pi(h'\oplus j')) =
  \kap(\pi(h\oplus j\oplus h'\oplus j')\oplus \pi(0)) =
  \kap(\pi(h\oplus h')\oplus\pi(j\oplus j')),
  \]
  so that $j\oplus j'\in A_{<k,<l}^\pi(i\oplus i',h\oplus h')$.

  The other relations follow in the same way, sometimes using
  (\ref{eq:7}) instead of~(\ref{eq:6}).
\end{proof}
The previous lemma has the following consequence for the sizes of the
sets $A_{\vartriangleleft k,\blacktriangleleft l}^\pi(i,h)$.
\begin{corollary}\label{cor:1}
  Let $\pi$ be a pseudo linear permutation, $i,h=0,\dots,2^n-1$, and
  $k,l=3,\dots,n$.
  \begin{enumerate}
  \item If $A_{<k,=l}^\pi(i,h)\neq\emptyset$ then
    $ |A_{<k,=l}^\pi(i,h)|=|A_{<k,<l}^\pi|.
    $\label{item:1}
  \item If $A_{=k,<l}^\pi\neq\emptyset$ then
    $ |A_{=k,<l}^\pi|=|A_{<k,<l}^\pi|.
    $
  \end{enumerate}
\end{corollary}
\begin{proof}
  Simply observe, that $A\oplus B\subseteq C$ and $B\neq\emptyset$
  imply that $|A|\leq|C|$ and use the relations proved in the previous
  lemma with the appropriate values for $i,i',h,h'$. E.~g.~choosing
  $i=i'$ and $h=h'$ in the first relation of Lemma~\ref{lem:8} gives
  $|A_{<k,=l}^\pi(i,\pi(h))| \leq |A_{<k,<l}^\pi|$, while $i=0$ and
  $h=0$ in the second relation yields $|A_{<k,<l}^\pi|\leq
  |A_{<k,=l}^\pi(i',\pi(h'))|$.  These two facts together
  prove~(\ref{item:1}).
\end{proof}
\begin{remark}
  Of course, there are many more of these relations, that can be
  proved in the same way. We will however only make use of these two
  relations here, so that we prove and formulate only the two.
\end{remark}

\begin{proposition}\label{prop:1}
  If $\pi$ is a pseudo linear permutation, then for fixed $i$, $h$ and
  $l$, the sets $A_k:=A_{<k,=l}^\pi(i,h)$ satisfy the following
  conditions. There exist numbers $p_k$ (possibly depending on $l$ but
  not on $i$ and $h$) satisfying $0\leq p_k\leq k$ and
  $p_{k+1}\in\{p_k,p_k+1\}$, and numbers $k_0$ (possibly depending on
  $i$, $h$ and $l$) such that $A_k$ is empty if and only if $k<k_0$
  while for $k\geq k_0$ we have $|A_k|=2^{p_k}$.
\end{proposition}
In other words, as soon as $k$ is so large that $A_k$ is non empty,
then either $|A_{k+1}|=|A_k|$ or $|A_{k+1}|=2|A_k|$. It is only the
set $A_{k_0}$ about which we have no information.
\begin{proof}
  The monotonicity in $k$ of the sets $A_k$ implies the existence of
  $k_0$. Now the second part of Corollary~\ref{cor:1} gives
  \[ |A_{<k,<l}^\pi|=|A_{<k-1,<l}^\pi| + |A_{=k-1,<l}^\pi|
  =2|A_{<k-1,<l}| \quad\mbox{if
    $A_{=k-1,<l}^\pi\neq\emptyset$.}
  \]
  Using the first part of Corollary~\ref{cor:1} we obtain that for
  $k>k_0$ the sets $A_k$ satisfy
  \[ |A_{k+1}|=
  \begin{cases}
    |A_k|  & \mbox{if $A_{=k-1,<l}^\pi=\emptyset$,}\\
    2|A_k| & \mbox{if $A_{=k-1,<l}^\pi\neq\emptyset$.}
  \end{cases}
  \]
  This proves the proposition.
\end{proof}
The next technical lemma prepares the proof of Theorem~\ref{thr:3}.
\begin{lemma}\label{lem:7}
  For $m\leq n$ let $q_m,\dots,q_n$ be numbers such that
  \[ q_k\geq0 \txt{and} q_{k+1}\in\{q_k,q_k+1\}.
  \]
  Then
  \[ \sum_{k=m}^n (-1)^k 2^{-q_k}  = (-1)^m \lambda,
  \]
  where $2^{-q_{m'}}/2\leq |\lambda|\leq 2\cdot2^{-q_{m'}}$ and $m'$ is
  the smallest of the numbers $m,\dots,n$ such that $m+m'$ is even and
  $q_{m'+1}=q_{m'}+1$. In particular $|\lambda|\leq 2$.
\end{lemma}
\begin{proof}
  Let
  \[ \mathbb{F}:=\{k<n:q_{k+1}=q_k+1,k+m \mbox{ even}\}
  \]
  and denote the elements of $\mathbb{F}$ by $k_1,\dots,k_\mu$.
  Obviously $m'=k_1$ and it follows by induction that
  \[ q_{k_\nu}\geq q_{m'}+\nu-1
  \txt{for $\nu=1,\dots,\mu$.}
  \]

  We now get for the sum in question
  \begin{align*}
    \sum_{k=m}^n (-1)^k 2^{-q_k} &=
    (-1)^m \sum_{\substack{k=m\\k+m \text{ even}}}
    (2^{-q_k} - 2^{-q_{k+1}}) + r,
  \end{align*}
  where
  \[ r=
  \begin{cases}
    0 & \text{if $n+m$ is odd,}\\
    2^{-q_n} & \text{if $n+m$ is even.}
  \end{cases}
  \]
  Using the definition of $\mathbb{F}$ and $q_{k_\nu}\geq
  q_{m'}+\nu-1$ for $\nu=1,\dots,\mu$ we obtain
  \[ \sum_{k\in\mathbb{F}} (2^{-q_k} - 2^{-q_{k+1}})
  = \frac12 \sum_{k\in\mathbb{F}}2^{-q_k}
  \begin{cases}
    \leq 2^{-q_{m'}} \sum_{\nu=1}^\mu 2^{-\nu} \leq 2^{-q_{m'}}, \\
    \geq \frac12 2^{-q_{k_1}} = \frac12 2^{-q_{m'}}.
  \end{cases}
  \]
  For the remaining part we obtain
  \[ \sum_{\substack{k=m\\k+m\text{ even}\\k\not\in\mathbb{F}}}
  (2^{-q_k} - 2^{-q_{k+1}}) =
  \sum_{\substack{k=m\\k+m\text{ even}\\k\not\in\mathbb{F}}}
  (2^{-q_k} - 2^{-q_k}) =0.
  \]
  This proves the assertion.
\end{proof}

If one is looking for a permutation $\pi$ maximizing $\alpha_n(\pi)$,
pseudo linear permutations might look like a good starting point,
since we have at least some information about the size of the sets
$A_{kl}^\pi(i,h)$. Given that the supremum over $h$ in the expression
for $\alpha_n(\iota)$ for the identical permutation $\iota$ is
actually attained if $h=i$, it might also look like a good idea to
first check the case $h=\pi(i)$. The next result tells us, that under
those assumptions we get again an upper estimate of only $\sqrt n$.
\begin{theorem}
  \label{thr:3}
  Let $\pi$ be a pseudo linear permutation. Then
  \[ \frac1{2^n} \sum_{i=0}^{2^n-1} \Big| \sum_{k,l=3}^n \pi(i)_l
  (-2)^{-k} |A_{<k,=l}^\pi(i,\pi(i))| \Big| \prec \sqrt n.
  \]
\end{theorem}
\begin{proof}
  First of all for pseudo linear permutations, the size of the sets
  $A_{<k,=l}^\pi(i,h)$ actually only depends on $i\oplus\pi^{-1}(h)$,
  since for $k,l\geq 3$
  \begin{align*}
    A_{<k,=l}^\pi(i,h)
    &= \{j:\kap(i\oplus j)<k, \kap(h\oplus\pi(j))=l\}\\
    &=\{j\oplus i: \kap(i\oplus j\oplus i)<k, \kap(h\oplus\pi(j\oplus
    i))=l\}.
  \end{align*}
  Now by pseudo linearity $\kap(h\oplus\pi(j\oplus
  i))=\kap(\pi(i\oplus\pi^{-1}(h))\oplus\pi(j))$. Hence
  \[ |A_{<k,=l}^\pi(i,h)|=|\{j:\kap(j)=k,
  \kap(\pi(i\oplus\pi^{-1}(h))\oplus\pi(j)) =l\}|.
  \]
  Therefore, the size of $A_{<k,=l}^\pi(i,\pi(i))$ is actually
  independent of $i$ and $h$ and hence equal to $|A_{<k,=l}^\pi|$.

  Next, we observe that
  \[ \bigcup_{l=3}^n A_{<k,=l}^\pi(i,h) =
  A_{<k,<n+1}^\pi(i,h)\setminus A_{<k,<3}^\pi(i,h).
  \]
  But it follows from Lemma~\ref{lem:5} that
  $|A_{<k,<n+1}^\pi(i,h)|=2^{k-1}$ and we trivially have
  $|A_{<k,<3}^\pi(i,h)|\leq 4$. Therefore
  \[ \Big| \sum_{k,l=3}^n (-2)^{-k} |A_{<k,=l}^\pi(i,h)| \Big| \leq
  \Big| \sum_{k=3}^n (-2)^{-k} 2^{k-1} \Big| + \sum_{k=3}^n
  2^{-k}\cdot 4 \leq \frac32.
  \]
  We can hence write
  \[
    \frac1{2^n}\! \sum_{i=0}^{2^n-1} \Big| \sum_{k,l=3}^n \pi(i)_l
    (-2)^{-k} |A_{<k,=l}^\pi(i,\pi(i))| \Big| \asymp
    \frac1{2^n}\!
    \sum_{i=0}^{2^n-1} \Big| \sum_{k,l=3}^n \big( i_l-\frac12\big)
    (-2)^{-k} |A_{<k,=l}^\pi| \Big|.
  \]
  It now follows from Proposition~\ref{prop:1}, that there are numbers
  $k_l$ and $p_{kl}$, such that $p_{k+1,l}\in\{p_{kl},p_{kl}+1\}$ and
  moreover $|A_{<k,=l}^\pi|=2^{p_{kl}}>0$ if and only if $k\geq k_l$.
  Therefore the sum in question becomes
  \begin{align*}
    \frac1{2^n}\! \sum_{i=0}^{2^n-1} \Big| \sum_{k,l=3}^n \pi(i)_l
    (-2)^{-k} |A_{<k,=l}^\pi(i,\pi(i))| \Big|
    &\asymp
    \frac1{2^n}\! \sum_{i=0}^{2^n-1} \Big| \sum_{l=3}^n
    \big(i_l-\frac12\big)
    \sum_{k\geq k_l} (-2)^{-k} 2^{p_{kl}} \Big|.
  \end{align*}
  Now, writing $t=i2^{-n}$, we get $i_l-1/2=r_l(t)/2$, where $r_l$
  denotes the $l$th Rademacher function and the sum becomes an
  integral over $t$:
  \[ \frac1{2^n}\sum_{i=0}^{2^n-1} \Big| \sum_{l=3}^n
  \big(i_l-\frac12\big)
  \sum_{k\geq k_l} (-2)^{-k} 2^{p_{kl}} \Big|
  = \frac12 \int_0^1 \Big| \sum_{l=3}^n
  r_l(t)
  \sum_{k\geq k_l} (-2)^{-k} 2^{p_{kl}} \Big|\, dt.
  \]
  By Khintchine's inequality this is equivalent to
  \[ \frac12 \int_0^1 \Big| \sum_{l=3}^n
  r_l(t)
  \sum_{k\geq k_l} (-2)^{-k} 2^{p_{kl}} \Big|\, dt \asymp
  \Big( \sum_{l=3}^n \Big| \sum_{k\geq k_l} (-2)^{-k} 2^{p_{kl}}
  \Big|^2 \Big)^{1/2}.
  \]
  Now the assertion follows using Lemma~\ref{lem:7} with $m=k_l$ and
  $q_k=k-p_{kl}$.
\end{proof}

\section{An optimizing strategy}
\label{sec:an-optim-strat}

Looking at the definition of $\INumdPerm_n$, its actual computation
can be done by a computer. However, as $n$ becomes larger, it quickly
turns out, that finding a maximum over $2^n!$ permutations is a task
to complex to be done in reasonable time. It is possible in this way
to find $\INumdPerm_1,\dots,\INumdPerm_4$ and the corresponding values
are listed in Table~\ref{tab:1}. As it turns out, in these cases there
is no improvement by taking permutations, that is, we have
\[ \INumdPerm_n=\alpha_n(\iota),\quad\mbox{for $n=1,2,3,4$.}
\]

\begin{table}
  \begin{center}
    \caption{$\INumdPerm_n$ for $n=1,2,3,4$ put in relation to
      $\alpha_n(\iota)$ and $\sqrt n$}
    \label{tab:1}
    \begin{tabular}{|r|l|l|l|l|}
      \hline
      $n$ &
      $\INumdPerm_n$ &
      $\INumdPerm_n/\sqrt n$ &
      $\alpha_n(\iota)$ &
      $\INumdPerm_n/\alpha_n(\iota)$ \\
      \hline
      $1$ & $0.25$           & $0.25$        & $0.25$        &
      $1.0$ \\
      $2$ & $0.5$            & $0.3535\dots$ & $0.5$         &
      $1.0$ \\
      $3$ & $0.5937\dots$    & $0.3428\dots$ & $0.5937\dots$ &
      $1.0$ \\
      $4$ & $0.6718\dots$    & $0.3359\dots$ & $0.6718\dots$ &
      $1.0$ \\
      \hline
    \end{tabular}
  \end{center}
\end{table}

So the question arises, whether one can, for larger $n$, find
permutations $\pi$ such that
\begin{equation}
  \alpha_n(\pi)>\alpha_n(\iota)=\sup_\pi\alpha_n(\pi).
  \label{eq:2}
\end{equation}
To see that this is indeed so, we now describe a general strategy to
find permutations $\pi$ such that $\alpha_n(\pi)$ becomes large. In
fact, I believe that using this strategy, one obtains permutations
$\pi$ such that
\[ \alpha_n(\pi)=\INumdPerm_n.
\]

Starting from an arbitrary permutation $\pi_0$ we want to find a
sequence of permutations $\pi_1,\pi_2,\dots$ such that
\[ \alpha_n(\pi_0)<\alpha_n(\pi_1)<\alpha_n(\pi_2)<\dots
\]
and $\pi_{k+1}$ is obtained from $\pi_k$ by a simple operation, such
as applying a transposition or a simple cycle. In fact it turns out
that cycles work best here and why this is so will be explained in
Proposition~\ref{prop:2}.

We can not proof that this strategy actually gives a value for
$\alpha_n(\pi)$ that comes even close to
$\INumdPerm_n=\sup_\pi\alpha_n(\pi)$, but there is some evidence that
it does.

For $h<i$ define the permutations $\gamma_{hi}$ and $\delta_{hi}$ by
\begin{align*}
  \gamma_{hi}&:
  \left( \begin{array}{rccrcl}
      0,\dots,h-1,&\hspace{3mm}h,\hspace{3mm}&h+1,& \hspace{-1mm}h+2,
      \dots,i-2, i-1,&           i,                 &i+1,\dots,2^n-1\\
      0,\dots,h-1,&            h,            &i,  & \hspace{-1mm}h+1,
      \dots,i-3, i-2,&           i-1,               &i+1,\dots,2^n-1
    \end{array} \right), \\
  \delta_{hi}&:
  \left( \begin{array}{rcrccr}
      0,\dots,h-1,&            h,            & h+1,              h+2,
      \dots,i-2,&i-1,&\hspace{2.5mm}i,\hspace{2.5mm}&i+1,\dots,2^n-1\\
      0,\dots,h-1,&            h+1,          & h+2,              h+3,
      \dots,i-1,& h, &              i,              &i+1,\dots,2^n-1
    \end{array} \right).
\end{align*}
In other words, $\gamma_{hi}$ is a cycle, taking $i$ to $h+1$, while
$\delta_{hi}$ is a cycle taking $h$ to $i-1$.

Consider now an arbitrary permutation $\pi$ of the set
$\{0,\dots,2^n-1\}$. Fix a number $i_0$ such that
$\pi(i_0)<\pi(i_0\oplus 1)$. (Since exchanging $i_0$ and $i_0\oplus1$
does not affect the value of $\alpha_n(\pi)$ at all, as one can easily
see, this is not really a restriction.) Then both permutations, $\pi'
:= \gamma_{\pi(i_0),\pi(i_0\oplus 1)} \circ \pi$ and
$\pi'':=\delta_{\pi(i_0),\pi(i_0\oplus 1)} \circ \pi$, are closer to
the identity than the original permutation $\pi$ in that for both the
images of $i_0$ and $i_0\oplus 1$ are consecutive numbers. The
strategy is now to pick $\pi_1$, among $\pi'$ and $\pi''$ such that
$\alpha_n(\pi_1)=\max(\alpha_n(\pi'),\alpha_n(\pi'')$. Hopefully
$\alpha_n(\pi_1)$ is also larger than $\alpha_n(\pi)$.

If this would always be so, we could simply repeat the process for the
next possible $i_0$ and end up with a permutation that always keeps
$i,i\oplus1$ together. We could then continue with the same argument
for groups of two numbers, $(i,i\oplus1)$ and $(i\oplus2,i\oplus3)$
and after $n$ passes we would arrive at the identity permutation,
maximizing $\alpha_n(\pi)$.

However, there are $i_0$ such that $\alpha_n(\pi_1) < \alpha_n(\pi)$.
For the sake of our strategy, we will then just keep fingers crossed,
leave the permutation as it was, continue with another value for
$i_0$, and see what happens.

Putting this strategy into a program, one finds local optimal values
$\INumdPermLoc_n$ for the function $\alpha_n(\pi)$ of which one can
hope, that they are at least close to the global maximum
$\INumdPerm_n$.

The results of the calculations are summarized in Table~\ref{tab:2}.
The final value for $\INumdPermLoc_n$ will depend on the randomly
chosen starting permutation. In the last two columns we have listed
the number of different starting permutations we used (the number of
runs for the program) and the number of starting permutations, for
which we obtained the specific maximal value $\INumdPermLoc_n$.

This shows, that the optimization gives indeed the maximum
$\INumdPerm_n$ for almost all starting permutations in the cases
$n=1,2,3,4$, where we can also compute the actual maximum. In the
other cases $n=5,6,7$ it seems that we can also find the maximum,
since for reasonably many starting permutations we end up with the
same value. Finally we include some calculations for the cases
$n=8,9$, which are not very reliable, since one run of the program for
the value $n=9$ for example already takes more than a day.

The results also show, that at least for the values of $n$ where our
information seems reliable, the asymptotic behavior of $\INumdPerm_n$
is rather close to $\sqrt n$ than to $n$.

\begin{table}
  \begin{center}
    \caption{Largest values $\INumdPermLoc_n$ we found for
      $\alpha_n(\pi)$}
    \label{tab:2}
    \begin{tabular}{|r|l|l|l|l|r|r|}
      \hline
      $n$ & $\INumdPermLoc_n$ &
      $\INumdPermLoc_n/\sqrt n$ & $\alpha_n(\iota)$ &
      $\INumdPermLoc_n/\alpha_n(\iota)$ &
      \#runs &
      \#successful \\
      \hline
      $1$ & $0.25$        & $0.25$        & $0.25$        &
      $1.0$         & $2000$ & $2000$ \\
      $2$ & $0.5$         & $0.3535\dots$ & $0.5$         &
      $1.0$         & $2000$ & $2000$ \\
      $3$ & $0.5937\dots$ & $0.3428\dots$ & $0.5937\dots$ &
      $1.0$         & $2000$ & $1789$ \\
      $4$ & $0.6718\dots$ & $0.3359\dots$ & $0.6718\dots$ &
      $1.0$         & $2000$ & $1260$ \\
      $5$ & $0.7509\dots$ & $0.3358\dots$ & $0.7421\dots$ &
      $1.0118\dots$ & $2000$ &  $392$ \\
      $6$ & $0.8203\dots$ & $0.3348\dots$ & $0.8046\dots$ &
      $1.0194\dots$ & $2000$ &  $124$ \\
      $7$ & $0.8847\dots$ & $0.3344\dots$ & $0.8632\dots$ &
      $1.0248\dots$ & $2000$ &  $4$ \\
      \hline
      $8$ & $0.9434\dots$ & $0.3335\dots$ & $0.9179\dots$ &
      $1.0277\dots$ & $750$ &  $2$ \\
      $9$ & $0.9970\dots$ & $0.3323\dots$ & $0.9697\dots$ &
      $1.0281\dots$ &  $8$ &  $1$ \\
      \hline
    \end{tabular}
  \end{center}
\end{table}

The C source code of the programs can be found in
Wenzel~\cite{wenzel03:_haar_banac} and can be downloaded from
\url{http://www.minet.uni-jena.de/~wenzel/habil/}.

We will now explain, why this strategy yields large values for
$\alpha_n(\pi)$. Denote by
\[ S^\pi(i) := \sup_h \Big| \sum_{j:\pi(j)\leq h} (-2)^{-\kap(i\oplus j)}
\Big|
\]
the supremum of the sums for fixed $i=0,\dots,2^n-1$ so that
\[ \alpha_n(\pi) = \frac1{2^n} \sum_{i=0}^{2^n-1} S^\pi(i).
\]
The next proposition says that the sum over `most' of the values for
$i$ becomes in fact larger when passing from $\pi$ to $\pi_1$, that is
the permutation such that
$\alpha_n(\pi_1)=\max\big(\alpha_n(\pi'),\alpha_n(\pi'')\big)$.
\begin{proposition}\label{prop:2}
  Let $i_0\in\{0,\dots,2^n-1\}$ be such that $\pi(i_0)<\pi(i_0\oplus 1)$
  and write $\pi':=\gamma_{\pi(i_0),\pi(i_0\oplus 1)} \circ\pi$ and
  $\pi'':=\delta_{\pi(i_0),\pi(i_0\oplus 1)} \circ\pi$. Then
  \[  \sum_{\substack{i=0\\i\not=i_0,i_0\oplus 1}}^{2^n-1} S^\pi(i)
  \leq
  \frac12 \Bigg(
  \sum_{\substack{i=0\\i\not=i_0,i_0\oplus1}}^{2^n-1} S^{\pi'}(i) +
  \sum_{\substack{i=0\\i\not=i_0,i_0\oplus1}}^{2^n-1} S^{\pi''}(i)
  \Bigg)
  \]
\end{proposition}
\begin{proof}
  Note first, that for $i\not=i_0,i_0\oplus1$ we have
  \begin{equation}\label{eq:18}
    \kap(i_0\oplus i) = \kap(i_0\oplus 1\oplus i).
  \end{equation}
  Given $i$, choose $h(i)$ such that
  \[ S^\pi(i) = \Big|
  \sum_{j:\pi(j)\leq h(i)} (-2)^{-\kap(i\oplus j)} \Big|.
  \]

  \begin{description}
  \item[Case 1 $h(i)<\pi(i_0)$ or $h(i)>\pi(i_0\oplus 1)$] Evidently
    for $\pi(j)\leq h(i)<\pi(i_0)$ we have
    \[ \pi'(j)=\pi(j)=\pi''(j)
    \]
    so that in the case $h(i)<\pi(i_0)$ we get
    \[ \{j:\pi(j)\leq h(i)\} = \{j:\pi'(j)\leq h(i)\} =
    \{j:\pi''(j)\leq h(i)\}.
    \]
    If $h(i)>\pi(i_0\oplus 1)$ then also
    \[ \{j:\pi(j)\leq h(i)\} = \{j:\pi'(j)\leq h(i)\} =
    \{j:\pi''(j)\leq h(i)\}
    \]
    since $\gamma_{\pi(i_0),\pi(i_0\oplus 1)}$ and
    $\delta_{\pi(i_0),\pi(i_0\oplus 1)}$ only move elements less
    than $h(i)$.

    In both cases, we obtain
    \begin{align*}
      S^\pi(i)
      &=
      \Big|
      \sum_{j:\pi(j)\leq h(i)} (-2)^{\kap(i\oplus j)}
      \Big|\\
      &=
      \frac12\Big|\sum_{j:\pi'(j)\leq h(i)} (-2)^{-\kap(i\oplus j)}\Big| +
      \frac12\Big|\sum_{i:\pi''(j)\leq h(i)} (-2)^{-\kap(i\oplus j)}
      \Big|\\
      &\leq
      \frac12\Big(S^{\pi'}(i) + S^{\pi''}(i)\Big).
    \end{align*}
  \item[Case 2 $\pi(i_0)\leq h(i)\leq\pi(i_0\oplus 1)$] It is clear
    from the definition of $\gamma$ that
    \[ \big\{j:\pi(j)\leq h(i) \big\} =
      \big\{j:\pi'(j)\leq h(i)+1\big\}
      \setminus
      \big\{i_0\oplus1\big\},
    \]
    while one can see from the definition of $\delta$ that
    \[ \big\{j:\pi(j)\leq h(i) \big\} =
      \big\{j:\pi''(j)\leq h(i)-1\big\}
      \cup
      \big\{i_0\big\}.
    \]
    Using~(\ref{eq:18}) we can write
    \begin{align*}
      \lefteqn{
        \sum_{j:\pi(j)\leq h(i)} (-2)^{-\kap(i\oplus j)}
        =
        \frac12 \Bigg(
        \sum_{j:\pi'(j)\leq h(i)+1} (-2)^{-\kap(i\oplus j)}
        - (-2)^{-\kap(i_0\oplus1\oplus j)} +} \hspace{2cm} \\
      &\hspace{2cm}
      \sum_{j:\pi''(j)\leq h(i)-1} (-2)^{-\kap(i\oplus j)}
      + (-2)^{-\kap(i_0\oplus j)} \bigg) \\
      &=
      \frac12 \Bigg(
      \sum_{j:\pi'(j)\leq h(i)+1} (-2)^{-\kap(i\oplus j)}+
      \sum_{j:\pi''(j)\leq h(i)-1} (-2)^{-\kap(i\oplus j)}
      \bigg)
    \end{align*}
    so by the triangle inequality also
    \begin{align*}
      S^\pi(i)
      &\leq
      \frac12\Big( S^{\pi'}(i)+S^{\pi''}(i)\Big).
    \end{align*}
  \end{description}
  Putting the two cases together and summing over the relevant $i$, we
  arrive at the assertion.
\end{proof}
As a consequence we see that, up to a perturbation, passing from $\pi$
to $\pi'$ or $\pi''$ indeed increases the value of $\alpha_n(\pi)$.
\begin{corollary}
  Writing
  \[ \delta(i):= S^\pi(i)-\frac{S^{\pi'}(i)+S^{\pi''}(i)}2,
  \]
  we have
  \[ \alpha_n(\pi) \leq \max\big(\alpha_n(\pi'),\alpha_n(\pi'')\big) +
  \delta(i_0) + \delta(i_0\oplus1).
  \]
\end{corollary}
Unfortunately we are not able to control the size of the perturbations
in a reasonable way so that the question about the actual asymptotic
behavior of the UMD constants of the summation operators remains open.


\end{document}